\ifx\shlhetal\undefinedcontrolsequence\let\shlhetal\relax\fi


\def\today{\ifcase\month\or January\or February\or March\or April\or
    May\or June\or July\or August\or September\or October\or November
    \or December\fi\space\number\day, \number\year}
\magnification=1200
\baselineskip=8mm
\font\bigbf=cmb10 scaled \magstep2
\def\nz{\hfil\break\noindent}
\def\nl{\nz}
\def\lz{\vskip6mm\noindent}
\def\hz{\vskip7.2pt\noindent}
\def\shelahtitle{\centerline{\bigbf How Special are Cohen and Random
Forcings}
\centerline{\bigbf i.e.
Boolean Algebras of the family of subsets}
\centerline{\bigbf of reals modulo meagre or
null}
\bigskip
Saharon
Shelah
\medskip
\centerline{\sl  Institute of Mathematics, The Hebrew University,
Jerusalem, Israel}
\centerline{\sl Department of Mathematics, Rutgers University,
New Brunswick, N.J., U.S.A.}}
\font\msam=msam10 scaled 1200
\font\msbm=msbm10 scaled 1200
\def\sp{{\rm sp}}
\def\lsp{{\rm \ell sp}}
\def\XXeq{\mathrel{\hbox{\msam\char"45}}}

\def\XX{\mathrel{\hbox{\msam\char"43}}}
\def\notXX{\;\>\rlap/\kern-5pt\XX}

\def\rk{\rangle}
\def\lk{\langle}

\def\lg{{\rm \ell g}}
\def\qed{\hfill\hbox{\msam\char'003}}
\def\Vdash{\kern2pt{\vrule depth 1pt height .8 em width .4pt}\kern
.1pt\vdash}
\def\notVdash{\;\>\rlap/\kern-5pt\Vdash}
\def\notvdash{\;\>\rlap/\kern-6pt\vdash}
\def\uhr{\hbox{\msam\char'026}}
\def\name#1{\lower 1em\rlap{\char'176}#1}
\def\lq{{\rm``}}
\def\rq{{\rm''}}
\def\and{{\, \&\,}}
\def\smallbox#1{\leavevmode\thinspace\hbox{\vrule\vtop{\vbox
   {\hrule\kern1pt\hbox{\vphantom{\tt/}\thinspace{\tt#1}\thinspace}}
   \kern1pt\hrule}\vrule}\thinspace}
\font\msbm=msbm10
\textfont9=\msbm
\font\msbms=msbm7
\scriptfont9=\msbms

\def\dbR{{\fam9 R}}
\def\bbr{\dbR}

\def\Min{{\rm Min\,}}
\font\fatone=cmssbx10
\def\bt{\hbox{\fatone t}}
\def\pr{{\rm pr\,}}
\def\conc{\widehat{~~}}
\def\Levy{{\rm Levy}}
\def\eqdf{\buildrel\rm def\over =}
\def\SP{{\cal P}}
\def\acad{{\rm partially \; supported\; by\; the\; basic\; research\;
fund,\; Israeli}\nl{\rm Academy}}
\def\Landau{{\rm partially\; sponsored\; by\; the\; Edmund\; Landau\;
Center\; for  \;
research}\nl{\rm in\; Mathematical
\; Analysis,\;supported\; by\; the\; Minerva\; Foundation\; (Germany)\;}.}
\today
\nl
 \shelahtitle
\bigskip
\footnote{}{dn 2/92, Publication 480, $\acad$, and $\Landau$}
The feeling that those two forcing notions
-Cohen and Random- (equivalently
the corresponding Boolean algebras
$\SP(\bbr)/$(meagre sets),
$\SP(\bbr)/$(null
sets)) are special,
was probably old and widespread.
 A reasonable interpretation is to show them unique,
or ``minimal'' or at least
characteristic in a
family of ``nice forcing'' like Borel.
We shall interpret ``nice'' as Souslin as suggested by
Judah Shelah [JdSh 292];
(discussed below).
We divide the family of Souslin forcing to two, and expect that:
among the first part, i.e. those
adding some non-dominated real,
Cohen is minimal
(=is below every one), while
among the rest random is quite characteristic even unique.
Concerning the second class we have weak results, concerning
the first class, our results look satisfactory.

Related is von Neumann's problem which in our language is:
\item{$(*)$} is there a
 ${}^\omega\omega$-bounding c.c.c.
forcing notion adding reals
which is not equivalent
to the measure algebra
(i.e. control measure problem)?

Velickovic (and , as I have  lately learnt, also Fremlin)
 suggests another problem (it says less on
forcings which are ${}^\omega\omega$-bounding but it says also much on
 the others).
\item{$(**)$} is there a c.c.c. forcing notion $P$ which adds new
reals and such that for every $f\in {}^{\omega}\omega\cap V^P$ there is
$h\in V$ such that $(\forall n)\vert h(n)\vert \le 2^n$, and $f(n)\in
h(n)$ for all $n$.
\nl
The version of it for Souslin forcing was our starting point.

We have two main results:
 one (1.14) says that Cohen forcing is ``minimal'' in the first class,
the other (1.10) says that all c.c.c.
 Souslin forcing have a property shared by Cohen forcing and Random real
forcing
(this is the answer to ($**$) for
Souslin forcing),
 so it gives a weak answer to the problem on how special is random
forcing,
but says much on all c.c.c. Souslin forcing.
 Earlier by  Gitik Shelah [Sh 412] , any
        $\sigma$-centered Souslin forcing notion add a Cohen real.
        We thank Andrzej Roslanowski for proof reading the paper very carefully
        correcting many and pointing out a flawed proof.

\bigskip \noindent
 {\bf \S1 A Souslin forcing which adds an unbounded real add a
Cohen real}
\bigskip\noindent
{\bf 1.1 Notation:}
0) $\lg(\eta)$ is the length of $\eta$.
\nl
1) $T$ denotes subtrees of ${}^{\omega>}\omega$,
i.e.,
$T\subseteq {}^{\omega>}\omega$ is non empty,
$[\nu\XX\eta\,\&\,\eta\in T\Rightarrow \nu\in T]$
and
[$\nu\in T\Rightarrow (\exists \eta\in T)(\nu\XX\eta)]$. For $\eta\in T$
let $T^{[\eta]}\eqdf\{\nu\in T:\nu\XXeq\eta$ or $\eta\XXeq\nu\}$ and let
$\lim
T=\{\eta\in {}^\omega\omega:\bigwedge_n\eta\uhr n\in T\}$
\nl
2) $\sp(T)=\{\eta\in T:
(\exists^{\ge 2} k)[\eta\conc\lk k\rk\in T]\}$,
\nl
 $\lsp (T)=\{\lg(\eta):\eta\in \sp (T)\}$.
\nl
3) $[A]^\mu=\{B\subseteq A:
|B|=\mu\}$,
$[A]^{<\mu}=\bigcup_{0\le \kappa<\mu}[A]^\kappa$
\nl
4) We say $T$ is $u$-large if:
$u\in [\omega]^{\aleph_0}$ and for some $n^*<\omega$: if
$n^*<n<m<\omega, n\in u, m\in u$ then
$[n, m)\cap \lsp(T)\not=\emptyset$.
\nl
5) We say $T$ is strongly $u-$large if:
$u\in [\omega]^{\aleph_0}$,
and for some
$n^*<\omega$,
if $n^*<n<m<\omega$,
$n\in u$,
$m\in u$ then
$(\forall \eta\in T\cap {}^n2)(\exists \nu)
[\eta\XXeq \nu\in \sp T\and  \lg \nu<m]$.
\nl
6) $O_k$ is a sequence of length $k$ of zeroes.
\nl
7) $(\forall^\infty n)$ means:
for every large enough $n<\omega$.
 $(\exists^\infty n)$ means for infinitely many $n<\omega$.
\nl
8) We say $T$ is $(u, \bar h)$-large if :
$u\in [\omega]^{\aleph_0}$,
$h_k:\omega\to \omega\setminus \{0, 1\}$,
$\bar h=\lk h_k:k<\omega\rk$
and for every
$k<\omega$,
 $T$ is
$(u, h_k)$-large which means:
for infinitely many
$n\in u$ we have:
\nl
 $n\le m\in u\,\and \, | u\cap m\setminus n|<h_k(n)\Rightarrow \Min
(\lsp(T)\setminus m)
<\Min(u\setminus(m+1))$.
\nl
Note

$(*)$ if $h_n=n$ for every $n<\omega$  this is equivalent to:
for every $k<\omega$, for some consequtive members 
$i_0<i_1<  \dots <i_k$  of $u$ , for every $\ell < k$
we have $ [i_\ell , i_{\ell+1} )  \cap \lsp(T)$ is not empty
\nl
9)  We say $\lk T_\ell:\ell<n\rk$ is $(u, \bar h)$-large if:
$u\in [\omega]^{\aleph_0},\bar h=\langle h_k:
k<\omega\rangle , h_k:\omega\to\omega\setminus\{0, 1\}$ and
for every $k<\omega$ for infinitely many $n\in u$ we have
\nl
$n\le m\in u\,\&\, |u\cap m\setminus n|
<h_k(n)\Rightarrow \bigwedge_\ell \Min(\lsp T_\ell\setminus m)<
\Min(u\setminus (n+1))$.
\nl
10) If $h_k=h$ for $k<\omega$ we write $h$ instead of $\bar h$.
\nl
11) We use forcing notions with the convention that larger means with
more information.
\nl
12) In a partial order (=forcing notion),
incompatible means have no common upper bound.
\lz
{\bf 1.2 Definition:}
A statement $\varphi(x)$ on reals is absolute if for every model $M$
extending $V$ with the same ordinals (mainly $M=V$ or a generic
extension) and
$N
$, a model of ZFC$^-$
(which is a transitive set or class 
of $M$
) with $\omega_1^M\subseteq N$ and
$a\in N$, we have
$N\models \varphi[a]$ iff $M\models \varphi[a]$.
\lz
{\bf 1.3 Definition:}
1) $P$ is a c.c.c. Souslin forcing notion if:
$P=(P, \le)$ is such that:
\nl
(a) there is a $\sum^1_1$-definition $\varphi^a$ of the set $P$
(which is $\subseteq \bbr$)
\nl
(b) there is a $\sum_1^1$-definition $\varphi^b$
of a partial order $\le$ on $P$
\nl
(c) there is a $\sum^1_1$-definition $\varphi^c$ of the relation $\lq
p, q$ incompatible in $P$'', see (1.1(12)) (hence it is $\Delta^1_1$,
as by the above it is $\Sigma^1_1$, now use Definition 1.3(1) (a)+(b),
it implies being compatible is $\sum^1_1$ hence being incompatible is
$\prod^1_1$).

\nl
(d) $(P, \le)$ satisfies the c.c.c..
\nl
2) Note: we do not distinguish strictly between $P$ and the three
$\sum_1^1$-formulas
$\varphi^a, \varphi^b, \varphi^c$ respectively
appearing in the definition.
\nl
3) $P$ is a Souslin forcing notion if $(a)+(b)$ holds.
\lz
{\bf 1.3A Remark:}
On (c.c.c.) Souslin forcing see Judah Shelah [JuSh 292] e.g.
\lz
{\bf 1.4 Claim:} 1) $\lq \varphi^a, \varphi^b, \varphi^c$ are
 $\sum^1_1$-formulas
 as in
$1.3(1)$'' is absolute.
\nl
2) For $P$ a c.c.c. Souslin forcing notion, $\lq\{r_n:n<\omega\}$ is a maximal antichain
of $P$''
is a conjunction of a $\sum^1_1$ and a $\prod^1_1$ statements.

\nl
3) Being a maximal antichain is absolute 
(even conjunction of $\prod^1_1$
and $\sum^1_1$)
hence so is ``being a
$P$-name of
 a member of ${}^\omega 2$ (or ${}^\omega\omega$)''.
\nl
4) If $P$ is a c.c.c.  Souslin forcing notion , $p_0 \in P $
and $P^* =^{df} \{ p \in P : P \models "p_0 \le p  \} $
(with the inherited order) then $P^*$  is a 
c.c.c.  Souslin forcing notion  too . 
\nl

{\it Proof:}
  E.g.
\nl
(2) The $\sum^1_1$ part is to say $\lq r_n\in P$'', so if $\lq x\in P
$'' is also $\prod^1_1$ then this statement is $\prod^1_1$; the
$\prod^1_1$ part is to say $(\forall x)\Big[x\notin P\vee
\bigvee_{n<\omega} (x, r_n$ compatible)\Big] (by Definition
1.3(1)((a)+(c)));
a third part is
$\bigwedge_{n<m<\omega}$
$(r_n, r_m$ incompatible) which are 
$\prod^1_1$ and $\Delta^1_1$
resp.

\nl
(3) Follows by (2)
\hfill\qed$_{1.4}$

\lz
{\bf 1.5 Lemma:}
 Assume $P$ is a c.c.c. Souslin forcing,
$\name r$ a $P$-name of a new member of ${}^\omega 2$.
 {\it Then} for some infinite $u\subseteq\omega$,
for every $p\in P$, the tree
$T_p[\name r]$
(see Definition 1.6 below) is $u$-large (see 1.1(5)).
\nl
\lz
{\bf 1.6 Definition:}
 $T_p[\name r]=\{\eta\in {}^{\omega>}2:
p\notVdash_P \lq \eta\not= \name r \uhr \lg \eta$''$\}$
(clearly it is a tree).
\nz
 Before we turn to proving Lemma 1.5, we prove:
\lz
{\bf 1.7 Claim:}
 1) For a given c.c.c. Souslin forcing notion $P$
(i.e. as in Definition 1.3(1)) and $P$-name $\name r$ of a
member of ${}^\omega 2$,
the conclusion of 1.5
 is an absolute statement
(actually $\sum^1_2$).
\nl
 2) The statement
 on $u,p$
(and also on $\name r$)
 that they is as required in 1.5, is  a $\sum^1_1$ statement.
\nl
3) Also ``$\name r$ is a $P$-name of a new real'' is absolute
in fact a $\prod^1_1$-statement ..
\nl
4) If $P$ is c.c.c. Souslin forcing notion,
above every $p\in P$ there are two incompatible conditions
{\it then} forcing with $P$ add a new real.
\hz
{\it Proof:}
1) Let $\name r$ be represented by
$\Big\lk\lk(p_i^\eta, \bt^\eta_i): \,i<\omega\rk$:
$\eta\in{}^{\omega>}2\Big\rk$
where
$\{p^\eta_i:i<\omega\}\subseteq P$
is a maximal antichain of $P$,
$\bt^\eta_i$ a truth value and $p_i^\eta\Vdash_P\lq\eta\XX\name r$ iff
$\bt^\eta_i$''.
For 1.5,
the failure of the statement
can be expressed by:
\smallskip

\item{$(*)$} $(\forall u)(\exists p)
\Big[u\subseteq \omega$ finite \ or \ $p\in  P \and
(\exists^\infty n\in u)(\forall \eta)[\eta\in {}^n2\and  
\eta\in T_p(\name r)
\Rightarrow \neg(\exists \nu)[\eta\XXeq \nu\and
\lg\nu<
\Min (u\setminus (n+1))\and
\nu\conc \lk 0\rk\in T_p(\name r)\and \nu\conc\lk 1\rk\in T_p(\name
r)]\Big]$.
\nl
Now the statement
``$\rho\in T_p(\name r)$''
is equivalent to
``$p\notVdash_P[\rho\notXX\name r]$''
which is equivalent to
\item{($**)$} $\bigvee_{i<\omega}(\bt^\rho_i=$ truth
$\and  p, p^\rho_i$ compatible.)

It is enough to show that $(*)$ is a $\prod^1_2$-statement hence it is
enough to show that inside the large parenthesis there is a
$\sum^1_1$-statement.  In $(*)$ inside the large parenthesis, ignoring
quantifications over $\omega$, we note that $\lq p\in P$'' is
$\sum^1_1$, and then we have to consider $(**)$, on which it is enough
to prove that it is a $\Delta^1_1$ statement [actually
we have three instances of it - all negatives].  By Definition
1.3(1)(c) it is $\prod^1_1$ and by Definition 1.3(1)(b) (and the
compatible meaning having a common upper bound) it is $\sum^1_1$.


2) The proof is included in the proofs of parts (1) and (3).

3) Easy. $\Big[$Why? the statement is
$(\forall p)[p\notin P \vee \bigvee_{\eta\in  {}^{\omega>} 2}$
 $[\eta\conc\langle 0\rangle \in T_p(\name r)\and\eta\conc \langle 1\rangle
\in T_p(\name r)] ]    $.
Now inside the 
parenthesis we have
$p \notin P$ which is $\prod^1_1$
 and  two instances of
$(**)$ which, as shown above, is a $\prod^1_1$-statement.$\Big]$

4) Easy,
e.g.  in $V^{Levy(\aleph_0 , 2^{\aleph_0} ) }$ we ask: is there $p
\in P$ such that $G_p =^{df} \{ q : q \in P^{V}, q \le p \} $
is a directed subset of $P^{V}$ , generic over $V$ , i.e.  not
disjoint to any maximal antichain of $P^{V}$ from $V$? By the
assumption if such $p$ exist , necessarily $G_p \notin V$,and by the
homogeneity we can find $Levy$-names $\name p , \name G_{\name p } $
of such objects
so in $V^{Levy} $ we can find a perfect set of such $G_p $'s , so the
$p$'s form an antichain of size continuum but this is absolute . So
there is no such $p$, letting $ \{ \{ p_{i,j} : i < \omega \} : j<
\omega \} $ list the maximal antichains of $P^V$ from $P^V$ (the list
in $V^{Levy}$ ), and we define a $p$-name $\name \eta \in {}^\omega
\omega$ : (in $V^{Levy}$): $\name \eta (n) = $the unique $m$ such that
$p_{n,m} \in G_P$, the generic subset of
$P^{V^{Levy}}$ . This is a $P$-name of a new real
(all in $V^{Levy(\aleph_0 , 2^{\aleph_0}) }  $ 
and by 
part (3) +1.4(2) its existence is absolute .
%
\hfill\qed$_{1.7}$
\lz
{\bf Remark:}
The use of $Q_D$ below can be replaced.
$Q_D$ is called Mathias forcing. See on it [Sh-b].
\lz
{\it 1.8 Proof of Lemma 1.5:}
Assume that the conclusion fails
(for $\name r $,
a $P$-name of a new member of ${}^\omega 2$,
which will be fixed until the end of the proof of Lemma 1.5).
For $D$ a filter on $\omega$
(containing the co-bounded subsets of $\omega$)
let $Q_D=\{(w, A):
w\subseteq \omega$ finite,
$A\in D$
and $\max(w)<\Min A$
(when $w\neq \emptyset)\}$
(and if $w\subseteq \omega$ is finite $A\subseteq \omega$
 we identify $(w, A)$ with $(w, A\cap(\max w, \omega)))$;
the order is defined by
 $(w_1, A_1)\le (w_2, A_2)$ iff
$w_1\subseteq w_2\subseteq w_1\cup A_1$,
$A_1\supseteq A_2$.
 Let $(w_1, A_1)\le_\pr(\omega_2, A_2)$
(pure extension) iff
$w_1=w_2$,
$A_1\supseteq A_2$.
Clearly $Q_D$ is a partial
order satisfying
the c.c.c. and     $\{q:q_0\le_\pr q\}$ is directed for each $q_0\in Q_D$.
Let $\name w=\cup\{w:(w, A)\in \name G_{{Q_D}}\}$,
clearly $\name w$ is
a $Q_D$-name and
any $G\subseteq Q_D$ generic over $V$
can be reconstructed from
$\name w[G]:G=\{(v, A)\in Q_D:v\subseteq \name w[G]\subseteq v\cup A\}$.
Without loss of generality CH holds (by claim 1.7(1),
e.g. force with $\Levy(\aleph_1, 2^{\aleph_0}$)),
hence we can choose $D$ as a Ramsey ultrafilter on $\omega$.
So as is well known that:
\nz
 $\otimes_1$ if $\name\ell<2$ is a $Q_D$-name
and $q\in Q$ {\it then} for some $q',
q\le_{\pr} q'\in Q_D$,
$q'$ forces a value to $\name\ell$.

So after forcing with $Q_D$,
the conclusion of 1.5 still fails (by claim 1.7(1)).
Hence for some $q^*\in Q_D$ and $Q_D$-names
$\name p, \name T$ , (remember that $\name r$
remains a $P$-name )
 we have $q^*$
$\Vdash_{Q_D}$ $\lq \name p\in P,\;\name r$ 
of ${}^\omega2$,
$\name T=T_{\name p}[\name r]$ is not
 $\name w$-large, such that: for arbitrarily large $n\in \name w$,
($\name w$-the $Q_D$-name) the
interval
 $\Big[n, \Min(\name w\setminus (n+1))\Big)$ is disjoint to
lsp$(\name T) $ " ; also we can assume
that
$\Vdash_{Q_D} \lq \name p \in P$
,$\name r$ 
remains a $P$-name of a new member 
of ${}^\omega 2$ and $\name T
=  T_{\name p}[\name r]  $ " .
  For $q\in Q_D$ let
$S[q]=:\{\eta\in {}^{\omega>} 2:$
for some $q', q\le_\pr q'$
and
$q'\Vdash_{Q_D} \lq \eta\in \name T$''$\}$;
note that $S[q]$ is also equal to
$\{\eta\in {}^{\omega>} 2:$
for no $q'$,
$q\le_\pr q'\in Q_D,
q'\Vdash_{Q_D} \lq \eta\notin \name T$''$\}$ (just apply $\otimes_1$).

Now note
\item{$\otimes_2$}
 $S[q]$ is a subtree (of ${}^{\omega>}2)$ and
if $q_1\le_\pr q_2$ (in $Q_D$)
{\it then} $S[q_1]\supseteq S[q_2]$, (in fact they are equal).
\nl
\item{$\otimes_3$} if
 $q^*\le q \in Q_D$ {\it then} for some
$q_1\ge q$ and $m$ we have:
 $S[q_1]$ has no splitting in any level $\ge m$.

\noindent
 Why?
let $n=\max w^q$;
  so $q$ forces that: for some $m$,
$m\in \name w$,
$m\ge n$,
and
$\Min[(\lsp T_{\name p}[\name r]\setminus m]\ge \Min [\name w\setminus
(m+1)]$.
 Before proving $\otimes_3$, repeatedly using $\otimes_1$ we can assume
\item{$\otimes_4$} if $m\in A^q$,
$v\subseteq m\cap A^q$,
$\eta\in {}^m 2$ then  the condition 
$(w^q\cup v\cup\{m\},
A\setminus (m+1)) \in Q_D$
forces
$(\Vdash_{Q_D})$
 a truth values to the following:

$(\alpha$) $\eta\in T_{\name p}(\name r)$

$(\beta)$ $(\exists \nu)\Big[\eta\XXeq \nu\in {}^{\omega>} 2
\and \lg \nu
<\Min(\name w\setminus(m+1)) \and
\nu\in \sp(T_{\name p}(\name r))\Big]$.

(Recall the  definition of  a Ramsey ultrafilter by game.)

By the sentence before the last,
for some
$m\in A^q$
and $v\subseteq A^q\cap m$
for every
$\eta^*\in {}^m2$,
if we get a positive answer for
$(\alpha)$ then we get a negative answer
for $(\beta)$;
let
$q'=(w^q\cup v\cup\{m\}$,
$A^q\setminus(m+1))$;
so $q'$ forces those two
statements.
Let for $k\in A^q\setminus (m+1)$,
$q^k=(w^q\cup v\cup\{m\},
A^q\setminus k)$
(so $q'\le_\pr q^k$)
and it forces
$(\Vdash_{Q_D})$
``every $\eta\in {}^m 2\cap T_{\name p}[\name r]$
has a unique extension in $T_{\name p}(\name r)\cap {}^k2$",
as required in $\otimes_3$.
\nl
 The rest of the argument will be used again so just note that
proving Claim 1.9 below enough for finishing the proof of 1.5.
\lz
{\bf 1.9 Claim:}
Assume $P$ is
a c.c.c. Souslin forcing,
$\name r$ a $P$-name of a new real,
$Q_D$, $S[q]$
(for $q\in Q_D$) chosen as above.
{\it Then} $\otimes_3$ above is impossible.
\hz
{\it Proof:} So assume $\otimes_3$ holds and we shall get,
eventually a contradiction.

For this end we define a forcing notion
$Q^*=Q^*_D$,
$Q^*_D=\{(\bar m, \bar q)$:
for some $n=
n(\bar q)=n(\bar m, \bar q)$
we have
 $\bar q=\lk q_\eta:\eta\in {}^n 2\rk$,
$q^*\le q_\eta\in Q_D$,
for each $\eta\in {}^n 2$ the sequence
$\lk \vert S[q_\eta]\cap {}^k 2\vert:
k<\omega\rk$
is bounded,
$\bar m=\lk m_\nu:\nu\in {}^{n>} 2\rk$,
$m_\nu<\omega$ and
{\it if}
$ \nu \conc \lk \ell\rk \XXeq \eta_\ell\in {}^n 2$ for
$\ell=0, 1$
 and  $k<\omega$ {\it then}
$\vert S[q_{\eta_0}]\cap S [q_{\eta_1}]\cap {}^k 2\vert\le m_\nu\}$.

The order is defined by
$(\bar m^1, \bar q^1)\le
(\bar m^2, \bar q^2)$
{\it iff}
$n(\bar q^1)\le n(\bar q^2)$,
$\bar m^1=\bar m^2\uhr
{}^{n(\bar q^1)>} 2$ and for
$\eta\in {}^{n(\bar q^2)} 2 $ we have
$q^1_{\eta\uhr n(\bar q^1)}\le
q^2_\eta$.

\nz
Clearly $Q^*_D$ satisfies the c.c.c. as for any
$(\bar m^*, \bar w^*) $ the
set
$\{(\bar m, \bar q)\in Q^*_D:
\bar m=\bar m^*,
w^{q_\eta}=w^*_\eta$ for every
$\eta\in {}^{n(\bar q)} 2\}$
 is directed (in $Q^*_D$).
\nz
  Also,
\item{$\otimes_5$} for every $n$,
$\{(\bar m, \bar q)\in Q^*_D$:
 $n(\bar q)\ge n\}$
is a dense (and open) subset of $Q^*_D$.
\nz
$\Big[$Why? it is enough to prove that for any given
$(\bar m^0, \bar q^0)\in Q^*_D$
with $n(0)\eqdf n(\bar q^0)$,
there is $(\bar m^1, \bar q^1)$ such that $(\bar m^0, \bar q^0)\le
(\bar m^1, \bar q^1)\in Q^*_D$ with
 $n(\bar q^1)=n(0)+1$;
let $q^1_\eta=q^0_{\eta\uhr n(0)}$,
and $m^1_\nu$ is: $m^0_\nu$ if
$\nu\in {}^{n(0)>} 2$ ,
and
$\max \{\vert S[q^0_\nu]\cap {}^{k}2\vert: k<\omega\}$
if $\nu\in {}^{n(0)}2$.
Check$\Big]$.

  For $G^*\subseteq Q^*_D$ generic over $V$,
 let for $\eta\in {}^{\omega>}2$,
$\name w_\eta[G^*]$ be
$\bigcup\{w^r:$
there is $(\bar m, \bar q)\in G^*$, $r=q_{\eta\uhr n}$ and
$n=n(\bar q)\le \lg (\eta)\}$ 
it is  well defined).
If $\eta\in ({}^\omega 2)^{V[G^*]}$ let
$\name w_\eta [G^*]$ be $\bigcup_{k<\omega} w_{\eta\uhr k} [G^*]$.

Also $Q^*_D$ adds a perfect set of generics
for $Q_D$, moreover:
in $V[G^*]$ for every $\eta\in ({}^\omega 2)^{V[G^*]},
\name w_\eta[G^*]$ defined above is generic
for $Q_D$ over $V$,
which means $G_\eta\eqdf\{(v, A):
v\subseteq \name w_\eta[G^*]\subseteq v\cup A\}$ is a generic subset
of $Q_D$ over $V$; this holds by $\otimes_6$ , $\otimes_7$
below .
\item{$\otimes_6$}
{\it if} $(\bar m, \bar q)\in Q^*_D$
and $\name \tau$ is
a $Q_D$-name of an ordinal {\it then}
we can find $ \bar q^1$
such that
$(\bar m, \bar q)\le (\bar m, \bar q^1)$,
$n(\bar q^1)=n(\bar q)$
and for every $\eta\in {}^{n(\bar q^1)}2$, the condition
 $q^1_\eta$ forces a value to $\name \tau$.
\nl
$\Big[$Why? let $\langle \eta_k:k<2^{n(\bar q)}\rangle$
list
${}^{n(\bar q)}2$.
 We now define
by induction on $k\le 2^{n(\bar q)}$,
a sequence
${\bar r}^k=\langle
r^k_\eta:\eta\in {}^{n(\bar q)} 2\rangle$
such that:
\nl
\phantom{such that:}
(a) $(\bar m, \bar r^k)\in Q^*_D$
\nl
\phantom{such that:}
(b) $(\bar m, \bar r^k)\le (\bar m, \bar r^{k+1})$
(i.e. $Q_D\models\lq r^k_\eta\le r^{k+1}_\eta$''
for $\eta\in {}^{n(\bar q)}2$)
\nl
\phantom{such that:}
(c) $\bar r^0=\bar q$
\nl
\phantom{such that:}
(d) $r^{k+1}_{\eta_k}$ forces a value
to $\name \tau$
(for the forcing notion $Q_D$).

If we succeed then
$\bar q^1\eqdf \bar r^{(2^{n(\bar q)})}$
is as required;
as
$\bar r^0$ is as required
the only problem is to find $\bar r^{k+1}$ being given
$\bar r^k$.
First we can find $m_k<\omega$,
such that no
$S[r^k_\eta]$ (for $\eta\in {}^{n(\bar q)}2)$ has a splitting node of level
$\ge m_k$,
and $m_k>n(\bar q)$.
Second we find
$r^{k, *}_{\eta_k}\in Q_D$ such that:
$Q_D\models$``$r^k_{\eta_k}\le_\pr r^{k, *} _{\eta_k}$''
and $r^{k, *}_{\eta_k}$ forces a truth value
to each statement of the form
 ``$\nu\in T_{\name p}[\name r]$''
for $\nu\in{}^{m_k\ge} 2$.
By the definition of $S[r^k_{\eta_k}]$ necessarily
$r^{k, *}_{\eta_k}\Vdash_{Q_D}$``$T_p[\name r]\cap^{m_k\ge}2\subseteq
S[r^k_{\eta_k}]$''.
Thirdly choose
$r^{k+1}_{\eta_k}\in Q_D$,
such that
$Q_D\models$
`` $ r^{k,*}_{\eta_k}\le r^{k+1}_{\eta_k}$''
and $r^{k+1}_{\eta_k}$ forces a value to $\name \tau$
(possible by density) and $S[r^{k+1}_{\eta_k}]$ has no splitting
above some level
(use $\otimes_3$).
Fourth,
let
$r^{k+1}_\eta=r^k_\eta$
for $\eta\in {}^{n(\bar q)} 2\setminus \{\eta_k\}$;
we still have to check
$\vert S[r^{k+1}_{\nu_0}]\cap S[r^{k+1}_{\nu_1}]
\cap {}^m2\vert \le m_\nu$
when $\nu\conc \langle \ell\rangle \XXeq \nu_\ell\in {}^{n(\bar q)} 2$;
by the induction hypothesis without loss of generality
$\eta_k\in \{\nu_0, \nu_1\}$,
so let $\eta_k=\nu_{\ell(*)}$.
If
$m\le m_k$ then
$S[r^{k+1}_{\nu_{\ell(*)}}]\cap{}^m 2\subseteq S[r^k_{\nu_{\ell(*)}}]$
and we are done
by the induction hypothesis.
If $m>m_k$,
by the choice of $m_k$,
$S[r^{k+1}_{\nu_{1-\ell(*)}}]$ has no splitting nodes of level
$\ge m_k$
hence
$\vert S[r^{k+1}_{\nu_{\ell(*)}}]\cap
S[r^{k+1}_{\nu_{1-\ell(*)}}]\cap {}^m 2\vert\le
\vert S[r^{k+1}_{\nu_{\ell(*)}}]\cap S[
r^{k+1}_{\nu_{1-\ell(*)}}]\cap
{}^{(m_k)} 2\vert$,
and use the previous sentence.
So we can carry the induction and
$\bar r^{(2^{n(\bar q)})}$ is as required in $\otimes_6\Big]$.
\nz
\item{$\otimes_7$}
{\it if} $(\bar m, \bar q)\in Q^*_D$
and $k < \omega $
{\it then}
we can find $ \bar q^1$
such that
$(\bar m, \bar q)\le (\bar m, \bar q^1)$,
$n(\bar q^1)=n(\bar q)$
and for every $\eta\in {}^{n(\bar q^1)}2$, the condition
 $q^1_\eta$ forces 
some $m_\eta$ to be in $w_\eta \setminus \{ k \} $
\nl
[Why? proof similar to that of $\otimes_6$ .  ] 

Now for every $\eta\in ({}^\omega 2)^{V[G^*]}$ we know that
$V[G_\eta]\models$
 ``$\name p[G_\eta]\in P$,
$\name r
$ is 
still
a $P$-name of a member of $(2^\omega )^{V[G_\eta][\name G_P]}$ and
$T_{\name p[G_\eta]}[\name r
]$ is not $\name w[G_\eta]$-large'',
by 1.7 this holds in $V[G^*]$ too.
A closer look shows that for $\eta\not= \nu$
(from $({}^\omega 2)^{V[G^*]}$) the tree
$T_{\name p[G_\eta]}[\name r]\cap
T_{\name p[G_\nu]}[\name r]$ has finitely many splittings.
So the conditions
 $\name p[G_\eta]$, $\name p[G_\nu]$ are incompatible in $P^{V[G^*]}$.
Contradiction to ``$P$ is c.c.c. Souslin and this is absolute (1.4(1))''.
\hfill\qed$_{1.9}$
\nl
\null\hfill\qed$_{1.5}$
\nz

Now we can answer Velickovic's question for Souslin forcings.
\lz
{\bf 1.10 Conclusion:}
Let $P$ be a c.c.c. Souslin forcing,
adding a new real.
\item{(1)}
The following is impossible: {\it for every} $P$-name of a new $\name
r\in {}^\omega\omega$ {\it for some} tree $T\subseteq
{}^{\omega>}\omega$, $T\in V$ and $p\in P$ we have $p\Vdash_P\lq \name r \in
\lim T$'', and $\bigwedge_{n<\omega}\vert T\cap {}^n \omega\vert \le 
2^n$, 
we can replace: ``for every $n$''   by  ``for infinitely many $n$''.
\item{(2)} The following is impossible:
{\it for some} $P$-name of a new $\name r\in {}^\omega\omega$ {\it for
every} strictly increasing $\{n_i:i<\omega\}\subseteq \omega$ from $V$
{\it for some} tree $T\subseteq {}^{\omega>}\omega$, $T\in V$ and $p\in P$ we
have $p\Vdash_P\lq \name r \in \lim T$'' and
$\bigwedge_{i<\omega}\vert T\cap {}^{n_i} \omega\vert \le
2^i$,
we can replace: ``for every $i$''   by  ``for inifinitely many $i$''.
\item{(3)}
The following is impossible: {\it for some} $P$-name $\name r $ of a
new member of ${}^\omega 2$
{\it for every} strictly increasing
$\{n_i:i<\omega\}\subseteq \omega$ from $V$ for some tree
$T\subseteq {}^{\omega>}2$ and $q$ we have 
$ q\in P $ and
$q\Vdash_P\lq \name r\in \lim T$'' and
$\bigwedge_{i<\omega}\vert T\cap {}^{n_i} 2\vert \le 
2^i$.
\item{(4)} The following is impossible:
{\it for some}
$r\in ({}^\omega 2)^{V^P}\setminus {}^\omega 2$
{\it for every}
strictly increasing sequence
$\lk n_i:i<\omega\rk\in V$ of
natural numbers,
for some tree $T\subseteq {}^{\omega>}2$ from $V$ we have:
$r\in \lim T$,
and
$(\exists^\infty i)\vert
T\cap
{}^{n_i}2\vert \le {}^i2$ or at least
$(\exists^\infty i)\vert T\cap {}^{(n_{i+1})}2\vert \le 2^{n_i}$.
\hz
{\it Proof:} (1), (2) follow by part (3) . Suppose that $\name r$ is a
$P$-name of a new member of ${}^\omega \omega$ .  Let $\name \eta_n$
be (the $P$-name ) $0_{\name r (n) +1 } \conc \lk 1 \rk $ and let
$\name r^*$ be the following $P$-name : the concatenation of $\name
\eta_0 , \name \eta_1 , \name \eta_2 , \dots $ .  By part (3) there is
a strictly increasing sequence $\lk n_i : i < \omega \rk $ of natural
numbers such that for no $q \in P$ and $T$ , does $q \Vdash \lq \name
r \in Lim (T) $ " and for inifinitely many $i < \omega$ we have
(2)  and replacing $\name r$  will give (1)  too .
\nl
(3) Follows from part (4).
\nl
(4) Let $\name r$ be a $P$-name of a new real.

By Lemma 1.5 for some infinite $u\subseteq \omega$
we have
\item{$(*)$} for every $p\in P$,
$T_p[\name r]$ is $u$-large
(see 1.1(4)).

\noindent
We now choose by induction on $i$, $n_i<\omega$, such that
$n_i>\sup\{n_j:j<i\}$ and $\vert (n_i, n_{i+1})\cap u\vert>
2^{n_i}+2$.  If $(4)$ fails for $\name r$ we apply the statement to
the sequence $\langle n_i:i<\omega\rangle$, so for some $p\in P$ and
subtree $T$ of ${}^{\omega>}2$ from $V$, we have:
\nl
\phantom{we have}
 (a) $p\Vdash_P$``$\name r\in \lim T$''
\nl
\phantom{we have}
(b) for infinitely many $i<\omega$, we have
$\vert T\cap^{(n_{i+1})}2\vert \le 2^{n_i}$.
\nl
By the choice of $u$, for some $j^*<\omega$, we know that $T_p[\name
r]$ has a splitting of level $\in [j_0, j_1)$ for each
$j_0,j_1\in u$, $j^*<j_0<j_1$.

So if $i>j^*$,
then $\vert T_p[\name r]\cap{}^{(n_{i+1})}2\vert$ is at least the number
of levels $<n_{i+1}$ of splitting nodes of
$T_p[\name r]$
which is $\ge \vert (n_i,n _{i+1})\cap u\vert $
which is $>2^{n_i}$.
But
$p\Vdash $``$\name r \in \lim T$''
implies $T_p[\name r]\subseteq T$
so $\vert T\cap {}^{(\eta_{i+1})}2\vert$
is $>2^{n_i}$,
( for every
$i<\omega$ such that
$i>j^*)$,
 this contradicts the choice of $T$ hence we finish.
\hfill\qed$_{1.10}$
\lz
{\bf 1.11 Remark:} This means that any c.c.c. Souslin forcing which is
${}^\omega\omega$-bounding is quite similar to the Random real forcing in
some sense.  More exactly every c.c.c. Souslin forcing has a property
shared by the Random real forcing and the Cohen forcing.
\lz
{\bf 1.12 Claim:}
1) Assume
\nl
\phantom{1.12 Cla}
 (a) $P$ is a forcing notion
\nl
\phantom{1.12 Cla}
(b) $\name r$ is a $P$-name of a member of ${}^\omega 2$.
\nl
\phantom{1.12 Cla}
(c) $\bar h=\lk h_n:n<\omega\rk,
h_n=n$ (i.e.
$h_n(i)=n$ for every $i<\omega$) and
$u\subseteq \omega$ is
infinite
\nl
\phantom{1.12 Cla}
(d) for every $p\in P$,
 for some $\eta\in T_p[\name r]$ the set
$\{k: \eta\conc O_{k-\lg\eta}\conc\lk 1\rk \in T_p [\name r]\}$
is
\nl
\phantom{1.12 Cla (d)}
 $(u, \bar h)$-large (see (*) of 1.1(8)).
\nl
{\it Then}
forcing with
$P$ add a Cohen real.
\nl
2) We can weaken (d) to
\nl
(d)$^-$ for every $p\in P$ for some
$n<\omega$,
 and
$\eta_0, \dots, \eta_{n-1}\in T_p[\name r]$ the set
 $\{k:$ for some $\ell<n$,
$\eta_\ell \conc O_{k-\lg \eta_\ell}\conc \lk 1\rk\in T_{p}[\name r
]\}$
is
$(u, \bar h)$-large.
\hz
{\it Proof:}
 1), 2) Let $u \setminus \{ 0 \} =\{n_i: 1\le i<\omega\}$,
$n_0\eqdf 0<n_1<n_2<\dots$,
let $\lk k(i, \ell):\ell<\omega\rk$
be such that
$i=\sum_\ell k(i, \ell) 2^\ell$,
$k(i, \ell)\in \{0, 1\},$ so
 $k(i, \ell)=0$ when
$2^\ell>i$.
Let $\rho^*_m=\lk k(i, \ell): \ell\le [\log_2 (i+1)]\rk$
where
$i=i_u(m)$ is the unique $i$ such that
$n_i\le m<n_{i+1}$.
We define a $P$-name $\name s$ (of a member of
$({}^\omega 2)^{V^P}):$
let $\{\name k_i:i<\omega\}$ list in increasing order
 $\{k<\omega:\name r(k)=1\}$
and
$\name s $ be $\rho^*_{k_0}\conc \rho^*_{{k_1}}\conc
\rho^*_{{k_2}}\conc\dots$.

Clearly by condition $(d)^-$,
for every $p\in P$
and $n<\omega$ we have
$p\notVdash$
``$\name r(k)=0$ for every $k\ge n$''.
Hence $\Vdash_P``\{k<\omega:
\name r(k)=1\}$
is infinite,
hence
$\Vdash_P\lq\name s \in {}^\omega 2$''.
It is enough to prove that $\Vdash_P$``$\name s $ is a Cohen real over
$V$''.
So let $T\in V$
be a given subtree of ${}^{\omega>}2$ which is nowhere dense,
i.e.
$(\forall\eta\in T)(\exists \nu)[\eta\XX \nu\in {}^{\omega>}2\setminus
T$],
and we should prove
$\Vdash_P$``$\name s\notin \lim T$''.
So assume
$p\in P$,
 $p\Vdash_{P}\lq\name s\in \lim T$''
and we shall get a contradiction.
Having our $p\in P$
we can apply $(d)^-$
(or ($d$) , which is stronger),
so we can find
$n<\omega$ and $\eta_{ 0 }, \dots, \eta_{n-1}\in T_p[\name r]$
as there
such that
$A=\{k<\omega$:
for some $\ell<n$,
$\eta_\ell\conc O_{k-\lg\eta_\ell}\conc
\langle 1\rangle \in
T_p[\name r]\}$
is
$(u, \bar h)$-large.

Let for each $\ell<n$,
$\{k_j^\ell: j<j_\ell\}$
list in increasing order
$\{k<\lg(\eta_\ell):
\eta_\ell(k)=1\}$ and let
$\rho^\ell=\rho^*_{k^\ell_0}\conc \rho^*_{k^\ell_1}\conc\dots
\rho^*_{k^\ell_{j_\ell-1}}$.
Now we can choose
by induction on $\ell\le n$,
a sequence $\nu_\ell\in {}^{\omega>}2$ such that:
$\nu_0=\langle \rangle$,
$\nu_\ell\XXeq \nu_{\ell+1}$ and
$\rho^\ell\conc \nu_{\ell+1}\notin T$
(each time use
``$T$ is nowhere dense'').

Next we choose
$m(*)\in A$
such that $\nu_n\XXeq \rho^*_{m(*)}$;
possible as $A$ is
$(n, \bar h)$-large
(check Definition 1.1(8):
the set $\{i_u(m):m\in A\}$
contains an interval
of length
$>2^{\lg(\nu_n)}$,
so by the definition of $\rho^*_m$,
some $m(*)$ in this interval is as required).
Now we can find $p_1\in P$
such that $p\le p_1$ and $p_1\Vdash_P$ ``for some
$\ell<n,$ $\eta_\ell\conc O_{m(*)-\lg\eta_\ell}\conc
\lk 1\rk\XXeq \name r$'' hence $p_1\Vdash_P$
``for some $\ell<n$, $\rho^\ell\conc
\rho^*_{m(*)}\XXeq  \name s$'',
so by the choice of $\nu_{\ell+1}$,
and as $\nu_{\ell + 1}  \XXeq \nu_n\XXeq \rho^*_{m(*)}$ we get
$p_1\Vdash_P$
``$\name s \notin \lim T$''
hence we get contradiction to:
$p\Vdash_P$``$\name s\in \lim T$'',
hence we finish proving $\Vdash_P$ ``$\name s$ is a Cohen real over
$V$.''
\hfill\qed$_{1.12}$
\lz
{\bf 1.13 Claim:}
Let $P$ be a c.c.c. Souslin forcing
\nl
1) $\lq P$ add a Cohen real'' is absolute (as well as ``$x$ is a $P$-name
of a Cohen real'').
\nl
2) $\lq x$ is a $P$-name of a dominating real'' is absolute.
\nl
3) $\lq P$ add a non dominated real'' is absolute
(as well as $\lq x$ is a $P$-name of a non dominated real'').
\nl
4)  for a given $\bar h$,
$\lq$ there is $u\in [\omega]^{\aleph_0}$ such that (d)
 of Claim 1.12(1) holds'' is absolute;
similarly $(d)^-$ of 1.12(2).
\nl
5)  $ \lq x$ is a $P$-name of a member of ${}^\omega \omega $ , 
dominating  $\eta_1 \in  {}^\omega \omega $ 
and not dominating 
$\eta_2 \notin  {}^\omega \omega $ "
is absolute (in fact , conjunction of
$\prod^1_1$  and  $\sum^1_1$ statements )
\hz
{\it Proof:}
1) Let $\varphi(x)$ say:
\itemitem{(a)}
 $x=\lk\lk p^\eta_i, \bt^\eta_i:i<\omega\rk:
\eta\in {}^{\omega>}\omega\rk$,
$p^\eta_i\in P$,
$\bt^\eta_i$ a truth value,
$\lk p^\eta_i:i<\omega\rk$
a maximal antichain,
( for each $\eta \in {}^{ \omega \ge } \omega$ )
\itemitem{(b)}
if $\eta,\, \nu\in {}^{\omega>}\omega$ and $p^\eta_i,
p^\nu_j$ are compatible
 {\it then}:
[$\eta\XXeq \nu \wedge
\bt_j^\nu$
truth
$\Rightarrow \bt^\eta_i$=truth] and
[$\eta, \nu $ are $\XX$-incomparable 
$\bigwedge  
\& \bt^\eta_i$=truth$\Rightarrow
\bt^\nu_j$=false].
\itemitem{(c)} for every $p\in P$ for some $\eta\in {}^{\omega>}\omega$
for every $\nu$,
$\eta\XXeq \nu \in {}^{\omega>}\omega$,
we have
$\bigvee_{i<\omega}(p, p^\nu_i$
compatible $\wedge$
$\bt^\nu_i$=truth).

Now by 1.4(1)+(2) part (a) is a conjunction of $\prod^1_1$ and $\sum^1_1$
 statements,
part (b)
is both $\prod^1_1$ and $\sum^1_1$ and part (c)
is $\prod^1_1$
(we use:
compatibility is both
$\prod^1_1$ and $\sum^1_1$
 and
($\forall p\in P)\;
[\dots]$ means $(\forall p)[p\notin P\vee \dots]$).
So $\varphi(x)$ is a conjunction of $\prod^1_1$ and
$\sum_1^1$ statements.
Now $\varphi(x)$ says
 ``$x$ represents a $P$-name of a Cohen real''
so $(\exists x)\varphi(x)$ which is a $\sum^1_2$ statement, express the
statement
``forcing with $P$ add a Cohen real.''
\nl
2) We repeat the proof of part (1) but clause (c) is replaced by:

(c)$^\prime$ {\it for every} $p\in P$
and $f\in ({}^\omega\omega)^V$
{\it there are} $q\in P$ and $n^*$ such that $p\le q$ and:

$(c)_{q, f, n^*}$ {\it if} $q, p^\eta_i$ are compatible,
$\bt^\eta_i$ truth and
$n^*\le n<\lg  (\eta ) $
{\it then} $f(n)\le \eta(n)$.
\nz
Now
 (c)$_{q, f, n^*}$ is a $\prod^1_1$ and $\sum^1_1$,
so $(c)'$ 
has the form
$(\forall p , f ) [ p \notin \vee   ( \exists q , n^* )
[ q \in P \& p \le q \& (c)_{q,f,n^*} ] $
which is 
$\prod^1_2$ hence ``$x$ is a $P$-name of a dominating real''
is an absolute statement.
\nl
3) Use the proof of Part (1) but
clause ($c$) is replaced by:

$(c)''$ for $p\in P$
for infinitely many $n$ the set
$\{\eta(n):\eta\in {}^{\omega>}\omega, \lg\eta>n,
i<\omega, \bt^n_i  = $ truth,
$p^\eta_i, p$ compatible$\}$ is infinite.
\nl
Now (c)$''$ is $\prod^1_1$ and we can finish as there.
\nl
4) The statement (d) and (d)$^-$ from 1.12  for given
$p, \name r, \bar h, u$ is a $\prod^1_1$ statement
(as by the proof of 1.7(1)  
$\nu \in T_p [ \name r ] $ is a $\prod^1_1$-statement
and a $\sum^1_1$-statement .)
\nl
5) Easier than the proof of (3)
\null\hfill\qed$_{1.14}$
\lz
{\bf 1.14 Conclusion:}
 If $P$ is a c.c.c. Souslin forcing notion
adding
$\name g\in {}^\omega\omega$
not dominated by any $f\in ({}^\omega\omega)^V$
{\it then} forcing with $P$
add a Cohen real.
\hz
{\it Proof:}
Without loss of generality, $\name g$ is strictly increasing.
Let $\name r=\{\name g(i):i<\omega\}$,
it is a subset of $\omega$ identified with its characteristic function.
We imitate the proof of Lemma 1.5 (using here 1.13(3) instead of 1.7
there) so as there without
 loss of generality there is a Ramsey ultrafilter
$D$ on $\omega$ and let the forcing notion $Q_D$ be as there.
Let $\bar h$ be as in Claim in 1.12 condition (c);
we ask:
\nl
$\otimes^1$
is there an infinite $u\subseteq \omega$ such that condition $(d)^-$ of claim
 1.12 holds?
\nl
If yes we are done
by claim 1.12 .
So from now on we asssume not.

Let  $G\subseteq Q_D$ be generic over $V$,
condition $(d)^-$ fails also  in  $V[G]$
(using absoluteness which holds by claim 1.13(4)),
in particular for $u=^{df} \name w [G]$.
For $p \in P^{V[G]}$  and  $\eta \in {}^{\omega >} 2$
we let

$C[  \eta , p ] =^{df}$
$\{ k : \eta \conc 0_{ k 
- \lg ( \eta ) }
\conc \lk 1 \rk  \in T_{p} [ {\name r } ]  \} $.


Hence for some $p^*\in P^{V[G]}$, 
(remeber  (*)  of 1.1(8)) :

$(*)_1 $ 
%
$V[G]\models \lq$ 
for every
 $j<\omega $ and $\eta^*_\ell \in T_{p^*}[\name r]$ for $\ell < j    $
for some $n^*$
  letting
$C
 =^{df}  $
$  \bigcup \{  C[ \eta^*_\ell , p^* ] : \ell < j \} $
 , for no $n^*+1$ consequtive members $i_0 < i_1<....<i_{n^*}$
of $u$ do we have  : for every $ m<n^*$
the sets  $[i_m , i_{m+1} ) $ , $C$ are not  disjoint \rq .





Let for $n<\omega$  , $  h_0 (n)$  be
like $n^*$ in $(*)_1$  for $\{ \eta^*_\ell : \ell < j \} =^{df}
T_{\name p^* } [ {\name r }  ] \cap {}^{n \ge } 2 $.
Let for  $n < \omega$  ,  $ h(n) =^{df}  \sum_{i \le n } (h_0 (i) + 1 )$.
Note that we  have:  $h$ is strictly increasing .

So  for some $q^* \in G $ (which is a subset of $Q_D$ ) , 
and $Q_D$-name
$\name p^*$ of a member of $P$ ,
we have:
$q^*$ forces that $\name p^*, 
{\name h}, { \name C } [ \eta , {\name p^*} ]$
(for $\eta \in {}^{\omega > } 2$ ,
)  are as above.

Wlog  $0 \in w^{q^*}$.

Let $ {   \name  w } = { \name w } \cup \{ 0 \} 
= \{  \name n_i : i < \omega \}$ be strictly increasing
, note $\name w$,$\name n_i $ are $Q_D$-names.

Wlog  

$\otimes^2$  for every  $k \in A^{q^* }$ and
  subset $v$ of $  \ A^{q^*} \cap k$
the condition  $(w^{q^*} \cup v \cup \{ k\} , A^{q^*}\setminus (k+1) )$
forces a value to the following:

(A) $T_{\name p^*} [ \name r ] \cap  {}^{k \ge} 2$,  say $t$

(B) truth value to $[\name n_i,\name n_{i+1})\cap C[ \eta , \name p^*
]=\emptyset$ for $\eta\in t$ and $i\le k+1$.



So

$\otimes^3$ Assume  $q_1 \ge q^*$, $ n < \omega$
and $q_1$ forces $ { \name h } (n) = n^* $  and 
$T_{ \name p^*} [ { \name r } ]  \cap  {}^{ n \ge } 2  = t$. If 
$ i_0<....<i_{n^*}$ are the first $n^*+1$ members of 
$A^{q_1}$, then for some $m<n^*$  the condition 
$q_2 =^{df} (w^{q_1} \cup  \{ i_0,...,i_{m-1} \} , A^{q_1}\setminus i_m)$
forces that: 

for no $\eta \in t $ and $j$ do we have:
$ i_{m-1} \le j  < Min ( {\name w }\setminus(i_{m-1} +1 ) ) $ and 
$ j \in  C[ \eta , { \name p^* } ]$.


Now we want to imitate the proof of 1.9 . We define a forcing
notion $Q^{**} = Q^{**}_D$  as follows :

a member of $Q^{**}_D$ has the form $( \bar m , \bar q )$
such that:

(1) $\bar q = \lk q_\nu :
\nu  \in {}^n 2 \rk $, $q^* \le q_\nu  \in Q_D$, $n=n(\bar q)$ ,

(2) for each $\nu   \in {}^n 2 $ the number
$i[q_\nu ] =^{df} Max ( w^{q_{\nu}})$ is well defined,

(3) $q_\nu $ forces (in $Q_D$) a value  $t_\nu $ to
$T_{\name p^*} [ {\name r } ]  \cap  {}^{n \ge } 2 $ and a value $h_\nu ( j )$
 to $ { \name h } (j) $ 
for $j \le n $  
and a value $ c_\nu ( \eta , \name p^* )$  to
$ {\name C } ( \eta , \name p^* ) \cap  i [ {q_\nu } ] $

(4) $q_\nu  \Vdash_{Q_D} \lq  $
if $\eta \in t_\nu $ then  for no $k$ do we have:
 $ i[q_\nu] \le  k  < 
Min ( \name w\setminus (i[  {q_\nu }]+1) )$ 
 and 
$ k \in  C[ \eta , \name p^* ] $ \rq


(5) $\bar m = \lk m_i: i\le n\rk $, $m_i<\omega$ and

(6) if 
$\nu \conc \lk j \rk \XXeq \nu_j \in {}^n 2 $ for  $j=0,1$ 
and  $\eta \in  t_{ \nu_0 } \cap t_{\nu_1}$
then 
$ c_{\nu_0} [ \eta  , p^* ]  \cap   c_{\nu_1} [ \eta , p^* ]$
is bounded by 
$\max\{m_{\lg  (\eta ) },m_{\lg (\nu ) } \}$

The order is defined  by : 
$(\bar m^1, \bar q^1)\le
(\bar m^2, \bar q^2)$
{\it iff}
$n(\bar q^1)\le n(\bar q^2)$,
$\bar m^1=\bar m^2\uhr
 ( n(\bar q^1) + 1  ) $ and for
$\eta\in {}^{n(\bar q^2)} 2 $ we have
$q^1_{\eta\uhr n(\bar q^1)}\le
q^2_\eta$.

We now continue as in the proof of 1.9 .


Why are the $ \name p^* [G_\nu ] $ for $\nu \in {}^\omega 2 $
pairwise incompatible?  If  $\nu_1,\nu_2 $ are
not equal , and ${\name p^*}[G_{\nu_1} ], {\name p^*}[G_{\nu_2}]$
are compatible in $P$
let $p$ be a common upper bound. We know that for some
$\eta  \in T_p [ {\name r } ] $ we have
$C[ \eta , p ] $ is infinite
as otherwise from $T_p [{\name r }]$ we can define
a function $f \in {}^\omega \omega$
such that 
$q \Vdash_P \lq {\name g } \le f $"
contradicting the assumption .

\null\hfill\qed$_{1.14}$
\nl
\centerline{$*$\qquad$*$\qquad$*$}

The question is whether a forcing adding half Cohen real (see below)
adds a Cohen real is due to Bartoszy\'nski and Fremlin, appears in [B].
\lz
{\bf 1.15 Definition:}
If $V\subseteq V^1$,
$r\in ({}^\omega\omega)^{V^1}$ we say that $r$ is a half 
Cohen real over $V$ {\it if} for every
$\eta\in ({}^\omega\omega)^V$,
for infinitely many $n<\omega$,
$r(n)=\eta(n)$.
\lz
{\bf 1.16 Conclusion:}
If a ccc Souslin forcing add a half Cohen real then it adds a Cohen real.
\hz
{\it Proof:} If $\name r$ is a $P$-name a member of ${}^\omega\omega$
which is (forced to be)
half Cohen over $V$,
then
$\Vdash_P$
 ``$\name r\in {}^\omega\omega$ is not dominated by any
``old''
$h\in {}^\omega\omega$ (i.e. $h\in ({}^\omega\omega)^V$''
trivially-you can
use the definition on $h+1$ to get strict inequality. Now use 1.14.
\lz

{\bf 1.17 Claim   : }
Assume $Q$  is a ccc Souslin forcing  and $\name r$
is a $Q$-name of a   member of ${}^\omega \omega$ .
Assume further that for some $\eta \in {}^\omega \omega$ , 
a Cohen real over $V$  and $G$ , a subset of
$Q^{V[\eta ] } $ generic over $V[ \eta ] $
we have 
 $\name r [G]$  dominate  $\eta$ 
i.e.  for every $n  <  \omega $  large enough 
$\eta (n) \le \name r [G] (n)$ .


Then   ( in $V$ )  $\name r $ is forced   by some  
$q \in Q^V$ to be a dominating real i.e. 
for every $G$  , a subset of $Q^V$ generic over $V$
to which $p $ belongs
and  $\rho \in  ( {}^\omega  \omega )^{V[ G]} $  
, for every $n < \omega $ large enough
$\rho (n) \le \name  r [G] (n)$

{\it Proof }: Assume that  not ;  then some $p \in  Q^{V[ \eta ] }$
force the   negation ,  i.e  $\name r $ is as above but 
the conclusion fail  .
By the homogeneity  of Cohen forcing
there is a Cohen name
$\name p $ of  such a condition  .

 Also in $V$ there is a maximal antichain 
$J$ of $Q$    and
sequence
 $\lk \rho_q  :  q \in J \rk$  such that
for each $q \in J$  either 
$q$ forces  that $\name r $ does not dominate
$\rho_q$   and  $\rho_q \in 
 ( {}^\omega  (\omega  \setminus \{ 0 \} )$
 OR  $q$ forces that  $\name r$ dominate every 
$\rho  \in   ({}^\omega \omega )^V$  and $\rho_q (n) = 0$
for every $n < \omega $
(all in $V$ ) . Now we can find $\rho^* \in   
{}^\omega   \omega $  dominating all $\rho_q$ for $q \in J$ . 

Also wlog  $\name p$ is above some some $q^* \in J$
so necessarily $\rho_{q^*} (n) \not=     0$ .

We define a forcing notion $R$ as follows: members has the form 
$(f,g)$  where for some $n=n(f) < \omega$, $f$ is a function 
from ${}^{n >} 2$  to  $\omega$ and $g$ is a function form 
$\omega$ to $\omega$ . The order is: $(f_1 , g_2 ) \le ( f_2 , g_2 )$
iff  
$n( f_1 ) \le  n( f_2 ) $ , 
$f_2$ extend $f_1$,  for every $n < \omega $ we have
$g_1 (n)  \le g_2 (n) $  and  for every  $ k $ satisfying
$ n( f_1 )  \le  k  < n( f_2 ) $  , for all except at most one
$ \nu \in {}^k 2 $ we have 
$ g_1 (  k ) \le f_2 ( \nu ) $ . 
\nl
{\bf  subclaim } Let $G$ be a  subset  of $R$  generic over $V$ . 
\item{(1)} In $V[G]$ , for every 
$\nu \in {}^\omega  2 $  , 
$\eta_\nu =^{df}  f \circ     \nu = \lk f( \nu \uhr    \ell ) :
\ell < \omega   \rk   $
is a Cohen real over 
$V$  , so 
\item{(2)} $p_{\eta_\nu }  =^{df} \name p [ \eta_\nu ] $  
is a member of 
$Q^{V [\eta ] } $  which by absoluteness is a subset of  $Q^{V[G]}$ .

\nl
Now  in $V[ \eta_\nu ] $  clearly $p_{\eta_\nu }$ 
forces that  $\name r$
dominate  $\eta_\nu$ but not $ \rho^* $ . By absoluteness this 
holds also  in $V[G]$ . Now if $\nu_1  ,  \nu_2$ are
distinct members of $({}^\omega  2 )^{V[G]}$   and 
$p_{\eta_{\nu_1}},
p_{\eta_{\nu_2}}$  are compatible in $Q^{V[G]}$ , let 
$p^* \in  Q^{V[G]}$
be a common upper bound  then it forces that $\name r$
dominates  $\rho^*$   , but it is also above a member $q^*$ 
of $J$  such that  $ \rho_{q^*} $ is not a sequence of zeroes.

By absoluteness we get a contradiction . 
\null\hfill\qed$_{1.17}$
\nl

\noindent {\bf References}

\item{[JdSh 292]}
H. Judah and S. Shelah, Souslin-forcing,
{\it J. Symb. Logic}, Vol. 53 no. 4 (1988):1188-1207.

\item{[Sh b]} Proper forcing, {\it Springer Lecture Notes},
940 (1982) 496 + xxix.

\item{[GiSh 412]} M. Gitik and S. Shelah, More on  simple forcing
notions, and forcings with ideals {\it Annals of Pure and
Applied Logic}.

\item{[B]} T. Bartoszy\'nski, Combinatorial aspects of measure
and category, {\it Fundamenta Mathematicae}, Vol. 127(1987):225--239.

\shlhetal

\bye